\documentclass[11pt]{amsart}
\title[$^{*}$-REGULARITY OF OPERATOR SPACE PROJECTIVE TENSOR PRODUCT ]{
$^{*}$-REGULARITY OF OPERATOR SPACE PROJECTIVE TENSOR PRODUCT OF C$^{*}$-ALGEBRAS
}
\usepackage[ansinew]{inputenc}
\usepackage{yfonts}
\usepackage{calligra}
\usepackage{amscd}
\usepackage{amssymb, latexsym}
\usepackage{mathrsfs}
\usepackage{fancyhdr}
\usepackage{enumerate}
\usepackage{textcomp, phonetic}
\usepackage{setspace}
\usepackage{latexsym}
\usepackage[all]{xy}
\usepackage{tikz}
\theoremstyle{plain}
\newtheorem{thm}{\sc Theorem}[section]
\newtheorem{cor}[thm]{\sc Corollary}
\newtheorem{lem}[thm]{\sc Lemma}

\newtheorem{prop}[thm]{\sc Proposition}

\newtheorem{defi}[thm]{\sc Definition}
\newenvironment{pf}{\noindent {\sc Proof:}}{\hfill $\Box$}
\begin{document}
\author[A. Kumar]{
AJAY KUMAR}
\address{Department of Mathematics\\
University of Delhi\\
Delhi\\
India.}
\email{akumar@maths.du.ac.in}
\thanks{}
\author[V. Rajpal]{
Vandana Rajpal}
\address{Department of Mathematics\\
University of Delhi\\
Delhi\\
India.}
\email{vandanarajpal.math@gmail.com}
\keywords{Operator space projective tensor norm, Enveloping $C^\ast$-algebra, $^{*}$-regularity.}
\subjclass[2010]{Primary 46L06, Secondary 46L07,47L25.}
\begin{abstract}
The Banach $^{*}$-algebra $A\widehat{\otimes}B$, the operator space projective tensor product of $C^{*}$-algebras $A$ and $B$, is shown to be $^{*}$-regular if Tomiyama's property ($F$) holds for $A\otimes_{\min}B$ and $A \otimes_{\min}B=A \otimes_{\max}B$, where $\otimes_{\min}$ and $\otimes_{\max}$
are the injective and projective $C^{*}$-cross norm, respectively.  However, $A\widehat{\otimes}B$ has a unique $C^{*}$-norm if and only if  $A\otimes B$ has. We also discuss the property ($F$) of $A\widehat{\otimes}B$ and $A\otimes_{h}B$, the Haagerup tensor product of $A$ and $B$.
\end{abstract}
\maketitle
\section{Introduction}
The concepts of $^{*}$-regularity and the uniqueness of $C^{*}$-norm have been extensively studied in Harmonic analysis for $L^{1}$-group algebras by J. Biodol~\cite{boi}, D. Poguntke~\cite{pogu}, Barnes~\cite{barne}, and others. Barnes in~\cite{Bar} studied these concepts in the context of $BG^{*}$-algebras. These results on tensor products were further improved by Hauenschild, Kaniuth, and Voigt~\cite{voi}.

Recall that if $A$ and $B$ are $C^{*}$-algebras and $u$ is an element in the algebraic tensor product $A\otimes B$, then the operator space projective tensor norm is defined to be\\
\hspace*{2.3 cm }$\|u\|_{\wedge}=\inf\{\|\alpha\|\|x\|\|y\|\|\beta\|:u=\alpha(x\otimes y)\beta\},$\\
where $\alpha\in M_{1,pq}$, $\beta\in M_{pq,1}$, $x\in M_{p}(A) $ and $y\in M_{q}(B)$, $p,q\in \mathbb{N}$, and
$x\otimes y=(x_{ij}\otimes y_{kl})_{(i,k),(j,l)}\in M_{pq}(A\otimes B)$.
The completion of $A\otimes B$ with respect to this norm is called the operator space
projective tensor product of $A$ and $B$ and is denoted by $A\widehat{\otimes}B$. It is well known that $A\widehat{\otimes}B$ is a Banach $^{*}$-algebra under the natural involution~\cite{kumar} and is a $C^{*}$-algebra if and only if $A$ or $B$ is $\mathbb{C}$. One of the main results about $^{*}$-regularity obtained in~\cite{voi},~\cite{pal} was that the Banach space projective tensor product of $C^{*}$-algebras $A$ and $B$ is $^{*}$-regular if their algebraic tensor product has unique $C^{*}$-norm and $A \otimes_{\min}B$ has property $(F)$. In Section 2, we prove a similar result for the operator space projective tensor product. Our proof is much shorter than the one given for Banach space projective tensor product case.  We also show that  $A\widehat{\otimes}B$ has a unique $C^{*}$-norm if and only if  $A\otimes B$ has. Using these results, we obtain several  Banach $^{*}$-algebras which are not $^{*}$-regular, e.g., $C^{*}_{r}(F_{2})\widehat{\otimes} C^{*}_{r}(F_{2})$, $B(H)\widehat{\otimes}B(H)$, and $B(H)/K(H)\widehat{\otimes}B(H)/K(H)$, where $C^{*}_{r}(F_{2})$ is the $C^{*}$-algebra associated to the left regular representation of the free group $F_{2}$ on two generators and $H$ being an infinite-dimensional separable Hilbert space, whereas the Banach $^{*}$-algebras  $C^{*}(G_{1})\widehat{\otimes} C^{*}(G_{2})$, $G_{1}$ and $G_{2}$ are  locally compact groups and $G_{1}$ is amenable, $K(H)\widehat{\otimes} K(H)$, $B(H)\widehat{\otimes} K(H)$, and  $B(H)\widehat{\otimes} K(H)+K(H)\widehat{\otimes} B(H)$ all are $^{*}$-regular. Section 3 deals with the $^{*}$-regularity of $A\widehat{\otimes}A$ with the reverse involution.  Finally, we introduce the notion of property $(F)$ for $A\widehat{\otimes}B$, and prove that if the Banach $^{*}$-algebra $A\widehat{\otimes}B$ has spectral synthesis in the sense of~\cite{ranjan}, then it satisfies property $(F)$. Similar result is also proved for the Haagerup tensor product of $C^{*}$-algebras.

\section{$^{*}$-Regularity And Unique $C^{*}$-norm
}
Throughout this paper, all $^{*}$-representations of $^{*}$-algebras are assumed to be normed and for any $^{*}$-algebra $A$, $Id(A)$ denotes the space of all  two-sided closed ideals of $A$. We begin by recalling some facts about $G^{*}$-algebras from~\cite{pal}. A $^{*}$-algebra $A$ is called a $G^{*}$-algebra if, for every $a\in A$, $\gamma_{A}(a)$ defined by
\[\gamma_{A}(a):=\sup\{\|\pi(a)\|: \;\; \text{$\pi$ a $^{*}$-representation of $A$}\},\]
is finite. This $\gamma_{A}$ is the largest $C^{*}$-seminorm on $A$; and the reducing ideal $A_{R}$ of $A$ (or $^{*}$-radical) is defined as $A_{R}=\{a\in A: \gamma_{A}(a)=0\}$, denote by $C^{*}(A)$, the completion of $A/A_{R}$ in the $C^{*}$-norm induced by $\gamma_{A}$. $C^{*}(A)$ together with the natural mapping $\varphi: A \to C^{*}(A)$ $(a\to a+A_{R})$ is called the enveloping $C^{*}$-algebra of $A$.
If $A_{R}=\{0\}$, i.e. if the points of $A$ are separated by its $^{*}$-representations, we say that $A$ is $^{*}$-reduced(or $^{*}$-semisimple).
Clearly, every Banach $^{*}$-algebra is a $G^{*}$-algebra. Also $A_{R}$, in this case, is a norm closed $^{*}$-ideal of $A$ and the quotient Banach $^{*}$-algebra $A/A_{R}$ is automatically $^{*}$-reduced. Note that, if $A$ is a $C^{*}$-algebra with respect to the norm $\|\cdot\|$, then $\gamma_{A}$ coincides with $\|\cdot\|$.

For a $^{*}$-algebra $A$, let $Prim^{*}(A)$ denote the set of all primitive ideals of $A$, i.e. set of kernels of topologically irreducible $^{*}$-representations of $A$. For a non-empty subset $E$ of $Prim^{*}(A)$, kernel of $E$ is defined to be\\
\hspace*{4 cm} $k(E)=\bigcap\{P:P\in E\}$,\\
 and for any subset $J$ of $A$, hull of $J$ relative to $Prim^{*}(A)$ is defined to be\\
\hspace*{3.9 cm} $h^{*}(J)=\{P\in Prim^{*}(A): P\supseteq J\}$.\\
We endow $Prim^{*}(A)$ with the hull-kernel topology (hk-topology), that is, for each subset $E$ of $Prim^{*}(A)$, its closure is $\overline{E}=h^{*}k(E)$. If $A$ is a $C^{*}$-algebra then we usually write $Prim(A)$ instead of $Prim^{*}(A)$.

In a similar manner, one can define the hk-topology on $Prime(A)$, the space of all prime ideals of $A$. Also recall that a $^{*}$-representation $\pi$ of a $^{*}$-algebra $A$ is called factorial if $\pi(A)^{''}$ (i.e., von Neumann algebra generated by  $\pi(A)$) is a factor. The set of kernels of factorial $^{*}$-representations of $A$ is called the factorial ideal space of $A$ and is denoted by $Fac(A)$. It is well known that the  kernel of a factorial $^{*}$-representation of a $^{*}$-algebra $A$ is a (closed) prime ideal, so that one can introduce the hull-kernel topology on $Fac(A)$.
\begin{defi}\emph{(~\cite{pal})}
A G$^{*}$-algebra $A$ is said to be $^{*}$-regular if the continuous surjection $\breve{\varphi} : Prim (C^{*}(A)) \to Prim^{*}(A)$ $(P \to \varphi^{-1}(P))$ is a homeomorphism, where $\varphi :A \to C^{*}(A)$ is the $C^{*}$-enveloping map of $A$.
\end{defi}
Equivalently, from (~\cite{chi}, Proposition 1.3), a Banach $^{*}$-algebra $A$ is $^{*}$-regular if and only if for any two non-degenerate $^{*}$-representations $\pi$ and $\rho$ of $A$, the inclusion $\ker \pi\subseteq \ker \rho$ will imply that $ \|\rho(a)\|\leq \|\pi(a)\|$ for all $a\in A$.

For $C^{*}$-algebras $A$ and $B$, it is known that $\|x\|_{\max}\leq \|x\|_{\wedge}$ for all $x\in A\otimes B$. So, there is a contractive homomorphism $i:A\widehat{\otimes}B \to A\otimes_{\max}B$ such that $i(a\otimes b)=a\otimes b$, for all $a\in A$, $b\in B$. Let $q$ be the canonical embedding from $A\otimes_{\max}B$ onto $A\otimes_{\min}B$. Then, by ~\cite{ranj}, $q\circ i$ is the canonical injection map from $A\widehat{\otimes}B$ to $A\otimes_{\min}B$. In particular, $i$ is injective. Also, given a closed ideal $I$ in $A\widehat{\otimes} B$, $I_{\max}$=cl($i(I)$) is a closed ideal in $A\otimes_{\max}B$, where closure is taken with respect to max-norm. From(~\cite{kumar}, Theorem 5), for a closed ideal $I$ of $A$ and $J$ of $B$, $I\widehat{\otimes} J$ is a closed ideal of $A\widehat{\otimes}B$. Therefore, cl($i(I\widehat{\otimes} J)$) is a closed ideal in $A\otimes_{\max}B$. That is, $I\otimes_{\max} J$ is a closed ideal in $A\otimes_{\max}B$, known as the product ideal.

Now recall that $A \otimes_{\min}B$ satisfies Tomiyama's property $(F)$ if the family $\{\phi\otimes_{\min} \varphi:\phi\in P(A),\;\varphi\in P(B)\}$ separates all closed ideals of $A \otimes_{\min}B$, where $P(A)$ and $P(B)$ denote the set of all pure states of $A$ and $B$, respectively. Using Proposition III. 1.5.10~\cite{bal}, Lemma 4.1~\cite{tak}, Proposition 3.6~\cite{ranja}, Theorem II.9.2.1~\cite{bal}, and Lemma 2.5~\cite{voi}, one can easily prove the following proposition  on the similar lines as given in Corollary 4.4~\cite{arc} for the Haagerup tensor product of $C^{*}$-algebras.
\begin{prop}\label{a12}
For $C^{*}$-algebras $A$ and $B$, the following conditions are equivalent.\\
\emph{(i)}$A \otimes_{\min}B=A \otimes_{\max}B$ and $A \otimes_{\min}B$ has property $(F)$. \\
\emph{(ii)}The mapping $K\to i^{-1}(K)$ is a homeomorphism from $Fac(A \otimes_{\max}B)$ onto $Fac(A \widehat{\otimes} B)$, where $i$ is the canonical map from $A\widehat{\otimes}B$ to $A\otimes_{\max}B$.\\
\emph{(iii)}The mapping $K\to i^{-1}(K) $ is a homeomorphism from $Prim(A \otimes_{\max}B)$ onto $Prim^{*}(A \widehat{\otimes} B)$, where $i$ is the canonical map from $A\widehat{\otimes}B$ to $A\otimes_{\max}B$.
\end{prop}
The following is an analogue of the result in~\cite{voi},~\cite{pal} for the operator space projective tensor product.
\begin{thm}\label{a113}
Let $A$ and $B$ be $C^{*}$-algebras, and suppose that $A \otimes_{\min}B$ has property $(F)$ and $A \otimes_{\min}B=A \otimes_{\max}B$. Then $A\widehat{\otimes}B$ is $^{*}$-regular. In particular, if $A$ or $B$ is nuclear then $A\widehat{\otimes}B$ is $^{*}$-regular.
\end{thm}
\begin{pf}
It can be seen easily that $i(A\widehat{\otimes}B)$ is $\|\cdot\|_{\max}$- dense in $A \otimes_{\max}B$, i.e., $\overline{i(A\widehat{\otimes}B)}=A \otimes_{\max}B$. Thus, by Theorem 10.1.11(c)~\cite{pal}, we get a unique $^{*}$-homomorphism $C^{*}(i):C^{*}(A\widehat{\otimes}B)\to C^{*}(A \otimes_{\max}B)= A \otimes_{\max}B $ that makes the following diagram commutative:
\\
\tikzset{node distance=3cm, auto}
\hspace*{3cm}
\begin{tikzpicture}
  \node (C) {$A\widehat{\otimes}B$};
  \node (P) [below of=C] {$C^{*}(A\widehat{\otimes}B)$};
  \node (Ai) [right of=P] {$A \otimes_{\max}B$};
  \draw[->] (C) to node {$i$} (Ai);
  \draw[->] (C) to node [swap] {$\varphi^{A\widehat{\otimes}B}$} (P);
  \draw[->] (P) to node [swap] {$C^{*}(i)$} (Ai);
\end{tikzpicture}\\
with $C^{*}(i)$ surjective. Also, by (~\cite{laur}, Theorem 4.8), $C^{*} (i)$ is an isometric isomorphism from $C^{*}(A\widehat{\otimes}B)$ onto $A\otimes_{\max}B$.\\
Since the canonical map $i:A\widehat{\otimes}B \to A \otimes_{\max}B$ is a $^{*}$-homomorphism, so Theorem 10.5.6~\cite{pal}, gives us a continuous map $\check{i}: Prim(A \otimes_{\max}B)\to Prim^{*}(A\widehat{\otimes}B)$
defined by $\breve{i}(P)=i^{-1}(P)$ for all $P\in Prim(A \otimes_{\max}B)$, and a commutative diagram:
\\
\tikzset{node distance=3.12cm, auto}
\hspace*{3cm}
\begin{tikzpicture}
  \node (C) {$Prim(A\otimes_{\max}B)$};
  \node (P) [below of=C] {$Prim(C^{*}(A\widehat{\otimes}B))$};
  \node (Ai) [right of=P] {$Prim^{*}(A\widehat{\otimes}B)$};
  \draw[->] (C) to node {$\breve{i}$} (Ai);
  \draw[->] (C) to node [swap] {$\breve{C^{*}(i)}$} (P);
  \draw[->] (P) to node [swap] {$\breve{\varphi}^{A\widehat{\otimes}B}$} (Ai);
\end{tikzpicture}
\\
i.e., $\breve{i}= \breve{\varphi}^{A\widehat{\otimes}B} \circ \breve{C^{*}(i)}$.

In order to show that $A\widehat{\otimes}B$ is $^{*}$-regular, suppose that $H$
is a closed subset of $Prim(C^{*}(A\widehat{\otimes}B))$. Since $\breve{C^{*}(i)}$ is a continuous map, so
$(\breve{C^{*}(i)})^{-1}(H)$ is closed in $Prim(A\otimes _{\max}B)$. By the given hypothesis and Proposition \ref{a12}, $\breve{i}$ is a homeomorphism from $Prim(A \otimes_{\max}B)$ onto $Prim^{*}(A \widehat{\otimes} B)$. So $\breve{i}((\breve{C^{*}(i)})^{-1}(H))$ is a closed subset of $Prim^{*}(A \widehat{\otimes} B)$.  We now claim that $\breve{C^{*}(i)}$ is a bijective map. It can be seen easily, using the surjectivity of $C^{*}(i)$, that $\breve{C^{*}(i)}$ is an injective map. To see the surjectivity, let $P\in Prim(C^{*}(A\widehat{\otimes}B))$, then $P=\ker \pi$, $\pi$ is an irreducible $^{*}$-representation of $C^{*}(A\widehat{\otimes}B)$. Let $\tilde{\pi}:=\pi\circ C^{*}(i)^{-1}$, then clearly $\ker \pi=\breve{C^{*}(i)}(\ker \tilde \pi))$ and $\tilde \pi$ is an irreducible $^{*}$-representation of $A\otimes_{\max}B$. So $\breve{i}= \breve{\varphi}^{A\widehat{\otimes}B} \circ \breve{C^{*}(i)}$ implies that $\breve{\varphi}^{A\widehat{\otimes}B}(H)$ is a closed subset of $Prim^{*}(A\widehat{\otimes}B)$; note that $\breve{\varphi}^{A\widehat{\otimes}B}$ is injective since $\breve{i}$
and $\breve{C^{*}(i)}$ both are bijective. Hence the result follows.
\end{pf}

For   an amenable locally compact group $G_{1}$  and infinite dimensional separable Hilbert space $H$, $C^{*}(G_{1})$ and $K(H)$  are nuclear, so $C^{*}(G_{1})\widehat{\otimes} C^{*}(G_{2})$, $B(H)\widehat{\otimes} K(H)$, $K(H)\widehat{\otimes} B(H)$, and $K(H)\widehat{\otimes} K(H)$, all are $^{*}$-regular Banach algebras, where $G_{2}$ is a locally compact group. Also, by (~\cite{pal}, Lemma 10.5.22), $B(H)\widehat{\otimes} K(H)+K(H)\widehat{\otimes} B(H)$ is $^{*}$-regular. Thus, every proper closed ideal of $B(H)\widehat{\otimes} B(H)$ is $^{*}$-regular.

A $G^{*}$-algebra $A$ has a unique $C^{*}$-norm if the Gelfand-Naimark norm $\gamma_{A/A_{R}}$ is the only $C^{*}$-norm that can be defined on $A/A_{R}$. Equivalently, from (~\cite{Bar}, Proposition 2.4), a reduced $BG^{*}$-algebra $A$ has a unique $C^{*}$-norm if and only if for every nonzero closed ideal $I$ of $C^{*}(A)$, $I\cap A $ is non-zero.

It is known, from~\cite{sinc}, that if $I$ is a non-zero closed ideal in $A\otimes_{\min}B$, then $I$ contains a non-zero elementary tensor. However, this may not be true for $A\otimes_{\max}B$ if $\|\cdot\|_{\min}\neq\|\cdot\|_{\max}$ on $A\otimes B$. This fact is essentially used in the following theorem.
\begin{thm}\label{a1322}
For $C^{*}$-algebras $A$ and $B$, $A\widehat{\otimes}B$ has a  unique $C^{*}$-norm if and only if $A\otimes B$ has.
\end{thm}
\begin{pf}
Suppose that $A\widehat{\otimes}B$ has a unique $C^{*}$-norm. In order to show the uniqueness of $C^{*}$-norm on $A\otimes B$, it is enough to show that $\|\cdot\|_{\min}=\|\cdot\|_{\max}$ on $A\otimes B$. Suppose, on the contrary,  $\|\cdot\|_{\min}\neq\|\cdot\|_{\max}$ on $A\otimes B$. Then  $\ker q$ is a non-zero closed ideal of $A\otimes_{\max}B$, where $q$ is the canonical $^{*}$-homomorphism from $A\otimes_{\max}B$ onto $A\otimes_{\min}B$. As in Theorem \ref{a113}, we obtain a unique isometric $^{*}$-isomorphism $C^{*}(i)$ from $C^{*}(A\widehat{\otimes}B)$  onto $A\otimes_{\max}B$. Therefore, the map $\phi:(Id(A \otimes_{\max}B), \tau_{\infty})\to (Id(C^{*}(A\widehat{\otimes}B), \tau_{\infty})$, given by $\phi(I)=C^{*} (i)^{-1}(I)$ for all $I\in Id(A \otimes_{\max}B)$, is a homeomorphism onto its image by (~\cite{beck}, Proposition 5), where $\tau_{\infty}$ is the topology on $Id(A)$~\cite{beck}. So $C^{*} (i)^{-1}(\ker q)$ is a non-zero closed ideal of $C^{*}(A\widehat{\otimes}B)$. Since $A\widehat{\otimes}B$ is $^{*}$-reduced, so $C^{*} (i)^{-1}(\ker q)\cap A\widehat{\otimes} B$ is a non-zero closed ideal of $A\widehat{\otimes}B$. Thus, by (~\cite{ranja}, Proposition 3.6), it would contain a non-zero elementary tensor, say $a\otimes b$, which further gives $C^{*} (i)(a \otimes b)\in \ker q$, i.e. $a \otimes b=0$, a contradiction. Hence $A\otimes B$ has a  unique $C^{*}$-norm.

Converse follows  by the same argument as that in (~\cite{pal}, Corollary 10.5.38).

\end{pf}

\begin{cor}\label{a132}
For $C^{*}$-algebras $A$ and $B$, if $A\widehat{\otimes}B$ is $^{*}$-regular, then $A\otimes B$ has a unique $C^{*}$-norm.
\end{cor}
\noindent From~\cite{voi}, $C^{*}_{r}(F_{2})\otimes C^{*}_{r}(F_{2})$ does not have unique $C^{*}$-norm, so $C^{*}_{r}(F_{2})\widehat{\otimes} C^{*}_{r}(F_{2})$ is not a $^{*}$-regular Banach algebra, where $C^{*}_{r}(F_{2})$ is the $C^{*}$-algebra associated to the left regular representation of the free group $F_{2}$ on two generators. Note that $C^{*}_{r}(F_{2})$ is non-nuclear simple $C^{*}$-algebra~\cite{bal}. Similarly, for any  $C^\ast$-algebra $A$ without the weak expectation property of Lance and  for a free group $F_{\infty}$ on an infinite set of generators, $C^\ast(F_{\infty})\widehat{\otimes}A$ is not $^{*}$-regular by Proposition 3.3~\cite{asym}. Also, by (~\cite{pisier}, Corollary 3.1), for an infinite dimensional separable Hilbert space $H$, $B(H) \otimes_{\min}B(H)\neq B(H) \otimes_{\max}B(H)$, so $B(H)\widehat{\otimes}B(H)$ is not $^{*}$-regular.
\begin{cor}
The maximal tensor product of two simple $C^\ast$-algebras need not be a simple $C^\ast$-algebra. Infact, if $A$ and $B$ are simple $C^\ast$-algebras such that $A\otimes_{\max}B$ is simple. Then $\|\cdot\|_{\min}=\|\cdot\|_{\max}$ on $A\otimes B$.
\end{cor}
\begin{thm}\label{a3}
Let $A$ and $B$ be  $^{*}$-reduced Banach $^{*}$-algebras, and let $\phi:A\to B$ be a bijective $^{*}$-homomorphism. Then $A$ is $^{*}$-regular if and only if   $B$ is.
\end{thm}
\begin{pf}
Suppose that  $A$ is a $^{*}$-regular Banach algebra. Since $\phi:A\to B$ is  an onto $^{*}$-homomorphism, so  Theorem 10.1.11(c)~\cite{pal} implies that there is a unique $^{*}$-homomorphism $C^{*}(\phi):C^{*}(A)\to C^{*}(B)$ that makes the following diagram commutative:
\[
\begin{CD}
A  @> \varphi^{A}>> C^{*}(A) \\
@V\phi VV    @VVC^{*}(\phi)V\\
B @>\varphi^{B}>> C^{*}(B)
\end{CD}
\]
with $C^{*}(\phi)$ surjective.
Also, by Theorem 10.5.6~\cite{pal}, we get the following commutative diagram:
\[
\begin{CD}
Prim^{*}(A)  @< \breve{\varphi}^{A}<< Prim(C^{*}(A)) \\
@A\breve{\phi} AA    @AA\breve{C^{*}(\phi)}A\\
Prim^{*}(B)  @<\breve{\varphi}^{B}<< Prim(C^{*}(B)).
\end{CD}
\]
Since  $A$ is $^{*}$-regular, $\breve{\varphi}^{A}$ is a homeomorphism. Now, we prove that both $\breve{\phi}$ and $\breve{C^{*}(\phi)}$  are bijective maps. Clearly,   the maps $\breve{C^{*}(\phi)}$ and
$\breve{\phi}$ are injective  since both $C^{*}(\phi)$ and $\phi$ were surjective. Before looking into the onto-ness of the maps $\breve{C^{*}(\phi)}$ and $\breve{\phi}$. We first show that $C^{*}(\phi)$ is an injective map. Since $\varphi^{B} \circ \phi=C^{*}(\phi)\circ \varphi^{A}$, so $\ker \phi=\ker C^{*}(\phi)\cap A$, which will be a non-zero ideal if  $\ker C^{*}(\phi)$ is  non-zero by the $^{*}$-regularity of $A$, a contradiction.  We now show that   $\breve{\phi}$ is an onto map. Let $P\in Prim^{*}(A)$ so $P=\ker \pi$, $\pi$ is an irreducible $^{*}$-representation of $A$. Define $\overline{\pi}:B\to B(H) $ as $\overline{\pi}(b)=\pi(\phi^{-1}(b))$, for all $b\in B$. Obviously, $\overline{\pi}$ is an irreducible $^{*}$-representation of $B$ and $\overline{\pi}\circ \phi=\pi$. Thus $\ker \pi=\ker \overline{\pi}\circ \phi=\breve{\phi}(\ker \overline{\pi})$, which shows that $\breve{\phi}$ is an onto map. Similarly, $\breve{C^{*}(\phi)}$  is an onto map.

In order to show the $^{*}$-regularity of $B$, let $H$ be a closed subset of $Prim(C^{*}(B))$. We  now  claim that $\breve{C^{*}(\phi)}(H)=h^{*}(\breve{C^{*}(\phi)}(k(H)))$.\\
Since $\breve{C^{*}(\phi)}$ is a bijective map, so\\
\hspace*{3 cm} $\breve{C^{*}(\phi)}(k(H))=\breve{C^{*}(\phi)}(\displaystyle\bigcap_{P\in H}P)$\\
\hspace*{5.384 cm}$=\displaystyle\bigcap_{P\in H} \breve{C^{*}(\phi)}(P) =k(\breve{C^{*}(\phi)}(H)).$\\
This proves \\
\hspace*{2.43 cm}$h^{*}( \breve{C^{*}(\phi)}(k(H)))=h^{*}(k(\breve{C^{*}(\phi)}(H)))$\\
\hspace*{5.42 cm}$\supseteq \breve{C^{*}(\phi)}(H)$.\\
Conversely, let $P\in h^{*}( \breve{C^{*}(\phi)}(k(H)))$ then $P\supseteq \breve{C^{*}(\phi)}(k(H))$. Since $\breve{C^{*}(\phi)}$ is a bijective map, so there exists a unique $Q\in Prim(C^{*}(B))$ such that $\breve{C^{*}(\phi)}(Q)=P$. Then, by using the bijectivity of $\breve{C^{*}(\phi)}$, we have $Q\in h^{*}(k(H))=H$ ($H$ is a closed set). Therefore $P\in \breve{C^{*}(\phi)}(H)$. This shows that  $\breve{C^{*}(\phi)}(H)$ is closed in $Prim(C^{*}(A))$. Since $\breve{\varphi}^{A}$ is homeomorphism and $\breve{\phi}$ is continuous, so $\breve{\phi}^{-1}(\breve{\varphi}^{A}(\breve{C^{*}(\phi)}(H)))$ is  closed  in $Prim^{*}(B)$. Thus $\breve{\varphi}^{B}(H)$ is a closed set in $Prim^{*}(B)$ since it equals $\breve{\phi}^{-1}(\breve{\varphi}^{A}(\breve{C^{*}(\phi)}(H)))$. Moreover, $\breve{\varphi}^{A}$, $\breve{C^{*}(\phi)}$, and $\breve{\phi}$ all are bijective, so $\breve{\varphi}^{B}$ is injective; hence $B$ is $^{*}$-regular.

Conversely, if $B$ is $^{*}$-regular, then so is $A$ by a similar argument, using $\phi^{-1}$.
\end{pf}

\begin{cor}\label{a121}
For an infinite dimensional separable Hilbert space $H$, \\$B(H)/K(H)
\widehat{\otimes}B(H)/K(H)$ is not a $^{*}$-regular Banach algebra.
\end{cor}
\begin{pf}
From ~\cite{ranjana}, we know that there exists an isometric isomorphism $\phi $ from $A:=(B(H)
\widehat{\otimes}B(H))/ (B(H)\widehat{\otimes}K(H)+K(H)\widehat{\otimes}B(H))$ to $B(H)/K(H)\widehat{\otimes}$ $B(H)/K(H)$, satisfying $\phi(x+(B(H)\widehat{\otimes}K(H)+K(H)\widehat{\otimes}B(H)))= q\widehat{\otimes}q(x)$, for all $x\in B(H)\widehat{\otimes}B(H)$, where $q$ is the quotient map from $B(H)$ to $B(H)/K(H)$. Clearly, this map $\phi$ is a bijective algebra $^{*}$-homomorphism. It is known from~\cite{ranjana} that  the  primitive ideals of $B(H)\widehat{\otimes}B(H)$ are $\{0\}$, $B(H)\widehat{\otimes}K(H)$, $K(H)\widehat{\otimes}B(H)$, and $B(H)\widehat{\otimes}K(H)+K(H)\widehat{\otimes}B(H)$. So $A$ has only one  primitive ideal  $(B(H)\widehat{\otimes}K(H)+K(H)\widehat{\otimes}B(H))/ (B(H)\widehat{\otimes}K(H)+K(H)\widehat{\otimes}B(H))$. Thus $A$ is $^{*}$-reduced. Now suppose that  $B(H)/K(H)\widehat{\otimes}B(H)/K(H)$ is $^{*}$-regular. Theorem \ref{a3} yields that $A$ is $^{*}$-regular and so is
$B(H)$ $\widehat{\otimes}B(H)$ by Theorem 10.5.15(d)~\cite{pal}, a contradiction.
\end{pf}
\section{Reverse Involution
}
Let $A$ be a $C^{*}$-algebra. On the Banach algebra $A \widehat{\otimes} A$, with the usual multiplication,
define the involution on an elementary tensor as $(a \otimes b)^{*} = b^{*}\otimes a^{*}$ for all $a, b \in A$.  This extends to $A \widehat{\otimes} A$, by the definition of operator space projective tensor norm and $A\widehat{\otimes} A$ becomes a Banach $^{*}$-algebra with this isometric involution, denoted by $A\widehat{\otimes}_{r} A$.
\begin{thm}\label{a211}
For a $C^{*}$-algebra $A$, $A\widehat{\otimes}_{r} A$ is $^{*}$-regular.
\end{thm}
\begin{pf}
Let $\pi$ and $\rho$ be non-degenerate $^{*}$-representations of $A\widehat{\otimes}_{r} A$, on the same Hilbert space $H$, with $\ker \pi\subseteq \ker \rho$. Suppose first that $A$ has an identity $1$. Define $\pi_{1}(a) := \pi(a \otimes 1)$ and $\pi_{2}(a) := \pi(1\otimes a)$, $a \in A$; clearly $\pi_{1}$ and $\pi_{2}$
are bounded representations from $A$ into $B(H)$ satisfying $\pi(a\otimes b)=\pi_{1}(a)\pi_{2}(b) = \pi_{2}(b)\pi_{1}(a)$ for all $a, b \in A $, and $\pi_{1}(a^{*}) = \pi_{2}(a)^{*}$ for all $a \in A$. Since every $^{*}$-representation of a Banach $^{*}$-algebra into a $C^{*}$-algebra is contractive and that $\widehat{\otimes}$ is a cross norm, so for a self-adjoint element $h\in A$, we get $ \|\exp (it\pi_{1}(h))\| = 1$ for all $t\in \mathbb{R}$. Thus $\pi_{1}(h)$ is a self-adjoint element of $B(H)$. Let $a\in A$, so $ a= h + ik$, where $h$ and $k$ are self-adjoint elements of $A$. Now, as in~\cite{AJAY}, $\pi_{1}(a^{*}) = \pi_{1}(a)^{*}$. This shows that $\pi_{1}$ is a $^{*}$-representation of $A$ and $\pi_{1}(a)^{*} = \pi_{2}(a)^{*}$ for all $a \in A$, and thus $\pi_{1}(a)=\pi_{2}(a)$ for all
$a \in A$. But $\pi(a\otimes b)=\pi_{1}(a)\pi_{2}(b) = \pi_{2}(b)\pi_{1}(a)$, so $\pi(a\otimes b)=\pi_{1}(ab)=\pi_{2}(ba)$, for all $a,b\in A$; similarly, $\pi_{2}$ is also a $^{*}$-representation of $A$. In a similar manner, we can define $^{*}$-representations $\rho_{1}$ and
$\rho_{2}$ of $A$ satisfying $\rho(a\otimes b)=\rho_{1}(a)\rho_{2}(b)=\rho_{2}(b)\rho_{1}(a)$ for all $a,b\in A$. Arguing as above, we have $\rho(a\otimes b)=\rho_{1}(ab)=\rho_{2}(ba)$, for all $a,b\in A$. Clearly, $\ker \pi_{1}\subseteq \ker \rho_{1}$.  As in ~\cite{tak}, it follows easily that  $\pi_{1}$, $\rho_{1}$, $\pi_{2}$, $\rho_{2}$, all are non-degenerate $^{*}$-representations of $A$. Therefore, by Proposition 1.3~\cite{chi}, we have \\
\hspace*{1.9 cm}$\|\rho(a\otimes b)\|=\|\rho_{1}(ab)\|\leq \|\pi_{1}(ab)\|=\|\pi(a\otimes b)\|$ for all $a,b\in A$.\\
Now for any $x=\displaystyle\sum_{i=1}^{n}a_{i}\otimes b_{i}$ in $A\otimes A$, clearly we have $\left\|\rho\left(\displaystyle\sum_{i=1}^{n}a_{i}\otimes b_{i}\right)\right\|\leq \left\|\pi\left(\displaystyle\sum_{i=1}^{n}a_{i}\otimes b_{i}\right)\right\|$. Since $A\otimes_{\wedge} A$ is dense in $A\widehat{\otimes}_{r} A$ and $^{*}$-representation from the Banach $^{*}$-algebra $A\widehat{\otimes}_{r} A$ to $B(H)$ is norm reducing, it follows easily that $\|\rho(x)\|\leq \|\pi(x)\|$ for all $x\in A\widehat{\otimes}_{r} A$. Hence if $A$ is a unital $C^{*}$-algebra, $A\widehat{\otimes}_{r} A$ is $^{*}$-regular by Proposition 1.3~\cite{chi}.

If $A$ does not have identity, consider the unitization $A_{e}$ of $A$. Clearly, $A$ is a closed ideal of $A_{e}$. Therefore $A\widehat{\otimes}_{r} A$ is a closed ideal of $A_{e}\widehat{\otimes}_{r} A_{e}$. Infact, it is a $^{*}$-ideal of $A_{e}\widehat{\otimes}_{r} A_{e}$. Thus $A_{e}\widehat{\otimes}_{r} A_{e}$ is $^{*}$-regular, so is $A\widehat{\otimes}_{r} A$ by Theorem 10.5.15~\cite{pal}. This completes the proof.
\end{pf}
\section{property $(F)$ for the operator space projective tensor product of $C^{*}$-algebras
}
Tomiyama in ~\cite{tomi} defined the concept of property $(F)$ for the minimal tensor product of $C^{*}$-algebras. Following ~\cite{tomi}, we define the property $(F)$ for the operator space projective tensor product of $C^{*}$-algebras and show that if $A\widehat{\otimes}B$, for any $C^{*}$-algebras $A$ and $B$, has spectral synthesis in the sense of~\cite{ranjan} then $A\widehat{\otimes}B$ satisfies property $(F)$. We also show that weak spectral synthesis and spectral synthesis in the sense of ~\cite{ranjan} coincides on $A\widehat{\otimes}B$.
\begin{defi}
Let $A$ and $B$ be $C^{*}$-algebras. We say that $A\widehat{\otimes}B$ satisfies property $(F)$ if the set of all product states $\phi\widehat{\otimes} \varphi$, where $\phi$, $\varphi$ are pure states of $A$ and $B$, respectively, separates the closed ideals of $A \widehat{\otimes}B$.
\end{defi}
From~\cite{ranj}, it is known that, for any $C^{*}$-algebras $A$ and $B$, the canonical map $i:A\widehat{\otimes}B \to A \otimes_{\min}B$ is an injective $^{*}$-homomorphism, so that we can regard $A\widehat{\otimes}B$
as a $^{*}$-subalgebra of $A \otimes_{\min}B$. Let $I$ be a closed ideal in $A\widehat{\otimes}B$ and $I_{\min}$ be the closure of $i(I)$ in $A \otimes_{\min}B$, in other words, $I_{\min}$ is the min-closure of
$I$ in $A \otimes_{\min}B$. Now associate two closed ideals with $I$ as:\\
$I_{l}$ = closure of span of all elementary tensors of $I$ in $A\widehat{\otimes}B$,\\
$I^{u}=I_{\min}\cap A\widehat{\otimes}B$, known as the lower and upper ideal associated with $I$, respectively. Clearly $I_{l}\subseteq I\subseteq I^{u}$. Following~\cite{ranjan}, we say that a closed ideal $I$ of $A\widehat{\otimes}B$ is spectral if $I_{l}= I= I^{u}$.

We start with the following  lemma which can be proved, using (~\cite{kumar}, Theorem 6),  on the same lines as done in Proposition 6.5~\cite{sinc} for the Haagerup tensor product of $C^{*}$-algebras.
\begin{lem}\label{a120}
Let $I$ and $J$ be closed ideals in $A\widehat{\otimes}B$ such that $I_{\min}=J_{\min}$. Then $J_{l}\subseteq I\subseteq J^{u}$ and $I_{l}\subseteq J\subseteq I^{u}$.
\end{lem}
We now relate the property $(F)$ of the operator space projective tensor product of $C^{*}$-algebras to the spectral synthesis of the closed sets of its primitive ideal space. For this, recall that for any Banach $^{*}$-algebra $A$, a closed subset $E$ of $Prim^{*}(A)$ is called spectral if $k(E)$ is the
only closed ideal in $A$ with hull equal to $E$, and we say that a Banach $^{*}$-algebra $A$ has spectral synthesis if every closed subset of $Prim^{*}(A)$ is spectral. From~\cite{ranjan},  for any $C^{*}$-algebras $A$ and $B$, the Banach $^{*}$-algebra $A\widehat{\otimes}B$ has spectral synthesis if and only if every closed ideal of $A\widehat{\otimes}B$ is spectral
\begin{thm}\label{a122}
Let $A$ and $B$ be $C^{*}$-algebras, and suppose that $A\widehat{\otimes}B$ has spectral synthesis. Then $A\widehat{\otimes}B$ satisfies property $(F)$.
\end{thm}
\begin{pf}
Let $I$ and $J$ be non-zero distinct closed ideals in $A\widehat{\otimes}B$. Then, by (~\cite{ranj}, Corollary 1), $I_{\min}$ and $J_{\min}$ are non-zero closed ideals in $A\otimes_{\min}B$. Suppose that $I_{\min}=J_{\min}$.  Using Lemma \ref{a120} and the fact that $A\widehat{\otimes}B$ has spectral synthesis, we get $I\subseteq J$, $J\subseteq I$, and thus $I=J$, a contradiction. This shows that $I_{\min}$ and $J_{\min}$ are non-zero distinct closed ideals in $A\otimes_{\min}B$. Now choose an irreducible $^{*}$-representation $\pi$ of $A\otimes_{\min}B$ such that $\pi(I_{\min})=0$ and $\pi(J_{\min})\neq 0$. Set $\tilde{\pi}:=\pi \circ i$, then by~\cite{ranjana}, $\tilde{\pi}$ is an irreducible $^{*}$-representation of $A\widehat{\otimes}B$ and clearly $\tilde{\pi}(I)=0$ and $\tilde{\pi}(J)\neq 0$. Let us denote the restriction of $\pi$ to $A$ and $B$ by $\pi_{1}$
and $\pi_{2}$, respectively; and $\tilde{\pi}_{1}$, $\tilde{\pi}_{2}$ be the restrictions of $\tilde{\pi}$ to $A$ and $B$.

Define a map $\theta_{\pi}:\pi(A \otimes B)\to \pi_{1}(A)\otimes \pi_{2}(B)$ as $\theta_{\pi}(\pi(a\otimes b))=\pi_{1}(a)\otimes \pi_{2}(b), \; a\in A,\;b\in B$. Then $\theta_{\pi}$ can be extended to a homomorphism $\tilde{\theta_{\pi}}$ from $\pi(A \otimes_{\min} B)$ onto $\pi_{1}(A) \otimes_{\min} \pi_{2}(B)$ and $\tilde{\theta_{\pi}}\circ \pi=\pi_{1}\otimes_{\min} \pi_{2}$ (see ~\cite{voi} and ~\cite{laur} for details). Note that $\tilde{\pi}(A \widehat{\otimes} B)\subseteq \pi(A \otimes_{\min} B)$. Since $\tilde{\pi}(J)\neq 0$, so choose $x \in J$ such that $\tilde{\pi}(x)\neq 0$. Suppose that $\tilde{\theta_{\pi}}(\tilde{\pi}(u))= 0$ for all $u\in J$. In particular, $\tilde{\theta_{\pi}}(\tilde{\pi}(x))= 0$, that is, $\pi_{1}\otimes_{\min} \pi_{2}(i(x))=0$. Since both $\pi_{1}\widehat{\otimes}\pi_{2}$ and $(\pi_{1}\otimes_{\min} \pi_{2})\circ i$ agree on $A\otimes B$, so by continuity $\pi_{1}\widehat{\otimes}\pi_{2}(x)=0$. Now we claim that $\tilde{\pi}_{1}=\pi_{1}$ and $\tilde{\pi}_{2}=\pi_{2}$. If $A$ and $B$ are unital, then $\tilde{\pi}_{1}(a)=\tilde{\pi}(a\otimes 1)=\pi(a\otimes 1)=\pi_{1}(a)$, for all $a\in A$, giving that $\tilde{\pi}_{1}=\pi_{1}$; similarly $\tilde{\pi}_{2}=\pi_{2}$. In the general case, if $\{e_{\lambda}\}$ and $\{f_{\mu}\}$ are the bounded approximate identities for $A$ and $B$, respectively, then for any $a \in A$, $\tilde{\pi}_{1}(a)=$s--$\lim\tilde{\pi}(a\otimes f_{\mu})=$s--$\lim\pi(a\otimes f_{\mu})=\pi_{1}(a)$ ~\cite{tak}, where s--lim denotes the strong limit. Thus $\tilde{\pi}_{1}=\pi_{1}$, similarly $\tilde{\pi}_{2}=\pi_{2}$, as claimed. By (~\cite{ranjana}, Theorem 7), $\ker \tilde{\pi}=A\widehat{\otimes}I_{2}+I_{1}\widehat{\otimes}B=\ker q_{I_{1}}\widehat{\otimes}q_{I_{2}}$, where $I_{1}=\ker \tilde{\pi_{1}}$, $I_{2}=\ker \tilde{\pi_{2}}$. By (~\cite{voi}, Lemma 2.1), $\ker \pi\subseteq \ker \pi_{1}\otimes_{\min} \pi_{2}$, giving that $\ker \pi\cap (A\widehat{\otimes}B) \subseteq \ker \pi_{1}\otimes_{\min} \pi_{2}\cap (A\widehat{\otimes}B)$, in other words, $\ker \tilde{\pi}\subseteq \ker \tilde{\pi}_{1}\widehat{\otimes} \tilde{\pi}_{2}$; that is, $\ker q_{I_{1}}\widehat{\otimes}q_{I_{2}} \subseteq \ker \tilde{\pi}_{1}\widehat{\otimes} \tilde{\pi}_{2}$. Suppose that the inclusion is strict. Let $K=\ker \tilde{\pi}_{1}\widehat{\otimes} \tilde{\pi}_{2}$, then $q_{I_{1}}\widehat{\otimes}q_{I_{2}}(K)$ is a non-zero closed ideal of $A/I_{1}\widehat{\otimes}B/I_{2}$ by (~\cite{ranjana}, Lemma 2).  So it must contain a non-zero elementary tensor, say $(a+I_{1})\otimes (b+I_{2})$(~\cite{ranja}, Proposition 3.6). Hence $a\otimes b\in K$, i.e., $\tilde{\pi}_{1}\widehat{\otimes} \tilde{\pi}_{2}(a\otimes b)=0$. So $\pi_{1}(a)\otimes \pi_{2}(b)=0$, i.e. either $\pi_{1}(a)=0$ or $\pi_{2}(b)=0$, a contradiction. Thus $\ker \tilde{\pi}_{1}\widehat{\otimes} \tilde{\pi}_{2}= \ker
q_{I_{1}}\widehat{\otimes}q_{I_{2}}=\ker \tilde{\pi}$. Therefore, $x\in \ker \tilde{\pi}$, which is not true. So $\tilde{\theta_{\pi}}(\tilde{\pi}(J))\neq 0$. Also note that $\tilde{\theta_{\pi}}(\tilde{\pi}(J))=
\pi_{1}\widehat{\otimes}\pi_{2}(J)$ and $\overline{\tilde{\theta_{\pi}}(\tilde{\pi}(J))}$, closure is taken with respect to min-norm, is a closed ideal in $\pi_{1}(A)\otimes_{\min}\pi_{2}(B)$.
Therefore, there exist $\phi\in P(\pi_{1}(A))$ and $\varphi\in P(\pi_{2}(B))$
such that $\phi\otimes_{\min}\varphi(\overline{\tilde{\theta_{\pi}}(\tilde{\pi}(J))})\neq 0$~\cite{tak},
so $\phi\otimes_{\min}\varphi(\tilde{\theta_{\pi}}(\tilde{\pi}(J)))\neq 0$, which further gives
$(\phi\circ \pi_{1} )\otimes_{\min}(\varphi \circ \pi_{2})(i(J))\neq 0$. Let $\sigma_{1}=\phi\circ \pi_{1}$ and $\sigma_{2}=\varphi\circ \pi_{2}$, then $\sigma_{1}\otimes_{\min}\sigma_{2}(i(J))\neq 0$, $\sigma_{1}\in P(A)$, $\sigma_{2}\in P(B)$. It is easy to see that both the maps $(\sigma_{1}\otimes_{\min}\sigma_{2}) \circ i$, $\sigma_{1}\widehat{\otimes}\sigma_{2}$ are continuous on $A\widehat{\otimes}B$ and agree on $A\otimes B$, giving that $\sigma_{1}\widehat{\otimes}\sigma_{2}(J)\neq 0$. Obviously $\sigma_{1}\widehat{\otimes}\sigma_{2}(I)=0$. Hence $A\widehat{\otimes}B$ has property ($F$).
\end{pf}

If $A$ or $B$ has finitely many closed ideals then $A\widehat{\otimes}B$ has spectral synthesis ~\cite{ranjan}. Thus, it satisfies property $(F)$. In particular, $B(H)\widehat{\otimes}B(H)$, $K(H)\widehat{\otimes}K(H)$ and $C_{0}(X)\widehat{\otimes}B(H)$ satisfy property $(F)$.

Recall that, for any $C^{*}$-algebras $A$ and $B$, an irreducible $^{*}$-representation of $A\otimes_{h}B$ is
a bounded homomorphism satisfying $\pi(a\otimes b)^{*} = \pi(a^{*} \otimes b^{*})$ for
all $a\in A$, $b \in B$, and for which $\pi(A\otimes_{h}B)$ is $\sigma$-weakly dense in
$B(H)$. From~\cite{sinc}, note that the irreduciblity of  $\pi$ is  equivalent to  $\pi(A\otimes_{h}B)'=\mathbb{C }I$.
Let $Prim(A\otimes_{h}B)$ denote the set of kernels of irreducible $^{*}$-representations
of $A\otimes_{h}B$, equipped with the hull-kernel topology. We shall  simply
write $h(I)$ for the hull of any ideal $I$ in $A\otimes_{h}B$, and we say  that $A\otimes_{h}B$ has spectral synthesis if every closed subset of  $Prim(A\otimes_{h}B)$ is spectral, where a closed subset $E$ of $Prim(A\otimes_{h}B)$ is called spectral if $k(E)$ is the only closed ideal in $A\otimes_{h}B$ with hull equal to $E$.

\begin{lem}
For $C^{*}$-algebras $A$ and $B$, if $A\otimes_{h}B$ has spectral synthesis then every closed ideal $I$ of $A\otimes_{h}B$ satisfies $I=I^{u}$.
\end{lem}
\begin{pf}
Let $I$ be a proper closed ideal of $A\otimes_{h}B$. As in (~\cite{ranjana}, Theorem 10), there exists an irreducible $^{*}$-representation $\sigma$ of $A\otimes_{h}B$ such that $I\subseteq \ker(\sigma)$, so that $\ker(\sigma)\in h(I)$, which shows that $h(I)$ is non-empty. Let $E=h(I)$, then $I=k(h(I))$ (see~\cite{ranjan} for details). Also $h(A\otimes_{h}B)=\emptyset$, and $k(\emptyset)=A\otimes_{h}B$, so that $A\otimes_{h}B=k(h(A\otimes_{h}B))$. This shows that every closed ideal $I$ of $A\otimes_{h}B$ is of the form $I=k(h(I))$, so $I=I^{u}$ by (~\cite{sinc}, Theorem 6.7).
\end{pf}

An argument similar to the one used to prove Theorem \ref{a122} establishes the following.
\begin{thm}
For  $C^{*}$-algebras $A$ and $B$, if $A\otimes_{h}B$ has spectral synthesis, then it satisfies property $(F)$.
\end{thm}
Recall that a closed subset $E$ of $Prime(A)$ is called weak spectral if $k(E)$ is the
only closed ideal in $A$ with hull equal to $E$. A Banach $^{*}$-algebra $A$ is said to
have weak spectral synthesis if every closed subset of $Prime(A)$ is weak spectral. The next result follows   on the similar lines as given  in ~\cite{ranjan}.
\begin{prop}\label{w1}
Let $A$ be a Banach $^{*}$-algebra having Wiener property. Then the following conditions are equivalent.\\
\emph{(i)} $A$ has weak spectral synthesis.\\
\emph{(ii)} Every $J\in Id(A)$ is of the form  $J = k(h^{p}(J))$, where $h^{p}(J)=\{P\in Prime(A) : P \supseteq J\}$.\\
\emph{(iii)} There is a one-one correspondence between the closed ideals of $A$ and the open subsets of $Prime(A)$.
\end{prop}

Let $Prim^{s}(A)$ denote the sobrification of the space  $Prim^{*}(A)$. Then, by (~\cite{spectral}, Lemma 1.1), $Prim^{s}(A)$ can be identified with the set of all semisimple prime ideals, where  $I \in Id(A)$ is said to be semisimple if it is an intersection of all primitive ideals of $A$ containing it.
\begin{prop}\label{w2}
Let $A$ be a Banach $^{*}$-algebra having Wiener property and $Prim^{s}(A)=Prime(A)$. Then $A$ has weak spectral synthesis if and only if $A$ has spectral synthesis.
\end{prop}
\begin{pf}
Suppose that $A$  has weak spectral synthesis. Then Proposition \ref{w1} yields that $Id(A)$ is isomorphic to the lattice of open subsets of $Prime(A)$, under the correspondence $I\to \{P\in Prime(A): P\nsupseteq I\}$. We know  that $Prim^{*}(A)$ and $Prim^{s}(A)$ have isomorphic lattices of open sets. So, by the given hypothesis,  $Id(A)$ is isomorphic to the lattice of open subsets of $Prim^{*}(A)$. Hence the result follows from (~\cite{ranjan}, Corollary 2.3).

Converse follows easily by the same argument using Proposition 4.6 above, Proposition 2.2 and Corollary 2.3~\cite{ranjan}.
\end{pf}

The following lemma can be proved  on the similar lines as done in Lemma 1.3~\cite{arc} for the Haagerup tensor product.
\begin{lem}\label{w3}
Let $I_{0}$ and $J_{0}$ be closed ideals in $C^{*}$-algebras $A$ and $B$, respectively. Let $S\subseteq Id(A)$, $T\subseteq Id(B)$ be such that $k(S)=I_{0}$ and $k(T)=J_{0}$. Then
$A\widehat{\otimes}J_{0}+I_{0}\widehat{\otimes}B=\bigcap\{A\widehat{\otimes}J+I\widehat{\otimes}B: I\in S, J\in T\}.$
\end{lem}
\begin{thm}
For $C^{*}$-algebras $A$ and $B$, weak spectral synthesis and spectral synthesis coincides on $A\widehat{\otimes}B$.
\end{thm}
\begin{pf}
Since $A\widehat{\otimes}B$ has the Wiener property (~\cite{ranjana}, Theorem 10), so by Proposition \ref{w2}, it suffices to show that  $Prim^{s}(A\widehat{\otimes}B)=Prime(A\widehat{\otimes}B)$. Let $P\in Prime(A\widehat{\otimes}B)$, so $P=A\widehat{\otimes} S +T \widehat{\otimes}B$, where $T\in Prime(A)$ and $S\in Prime(B)$ (~\cite{ranjana}, Theorem 6). Since  $T\in Prime(A)$ and $S\in Prime(B)$, so $T=k(h^{*}(T))$ and $S=k(h^{*}(S))$. Therefore, by Lemma \ref{w3}, $P=\bigcap \{A\widehat{\otimes}P_{2}+P_{1}\widehat{\otimes}B: P_{1}\in h^{*}(T), P_{2}\in h^{*}(S)\}$. Since the product ideals in $A\widehat{\otimes}B$ are closed (~\cite{kumar}, Theorem 5), so for $P_{1}\in h^{*}(T)$, $P_{2}\in h^{*}(S)$, $A\widehat{\otimes}P_{2}+P_{1}\widehat{\otimes}B\supseteq P$. Thus, every prime ideal of $A\widehat{\otimes}B$ is an intersection of all the primitive ideals of $A\widehat{\otimes}B$ containing it, so is semisimple.
 \end{pf}

\end{document}